\documentclass[reqno,12pt,letterpaper]{amsart}
\usepackage{overpic}
\usepackage{amsthm}
\usepackage{graphicx}
\usepackage{times}
\usepackage{geometry}
\usepackage{fullpage}
\usepackage{amsmath}
\usepackage{amssymb}
\usepackage{graphs}
\usepackage{varioref}
\usepackage{hyperref}
\usepackage{units}

\newtheorem{thm}{Theorem}

\newtheorem{Def}{Definition}
\newtheorem{proposition}{Proposition}
\newtheorem{claim}{Claim}
\newtheorem{lemma}{Lemma}
\newtheorem{conjecture}{Conjecture}
\newtheorem{case}{Case}

\newcommand{\R}{\mathbb{R}}

\newcommand{\FTC}{\operatorname{FTC}}

\newcommand{\Len}{\operatorname{Len}}

\title{A Class of Curves In Every Knot Type Where Chords of High Distortion are Common}

\author[Mullikin]{Chad A. S. Mullikin}
\address{Department of Mathematics, Spring Hill College,
Mobile, AL 36608}
\email{chadm@math.uga.edu}


\begin{document}

\begin{abstract} 
  The distortion of a curve is the supremum, taken over distinct pairs of points of the curve, of the ratio of arclength to spatial distance between the points. Gromov asked in 1981 whether a curve in every knot type can be constructed with distortion less than a universal constant C. Answering Gromov's question seems to require the construction of lower bounds on the distortion of knots in terms of some topological invariant. We attempt to make such bounds easier to construct by showing that pairs of points with high distortion are very common on curves of minimum length in the set of curves in a given knot type with distortion bounded above and distortion thickness bounded below.
\end{abstract}

\maketitle    % Creates title page, copyright page if any, and approval page.

\section{A New Point Of View On Gromov's Distortion}

\begin{figure}[htbp]
  \begin{tabular}{cc}
    \hspace{0.0in}\begin{overpic}[width=1.5in]{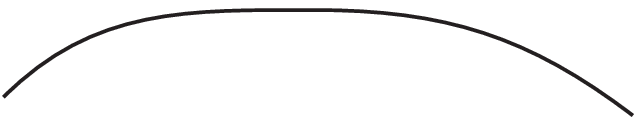}\end{overpic}&
    \hspace{0.5in}\begin{overpic}[width=1.5in]{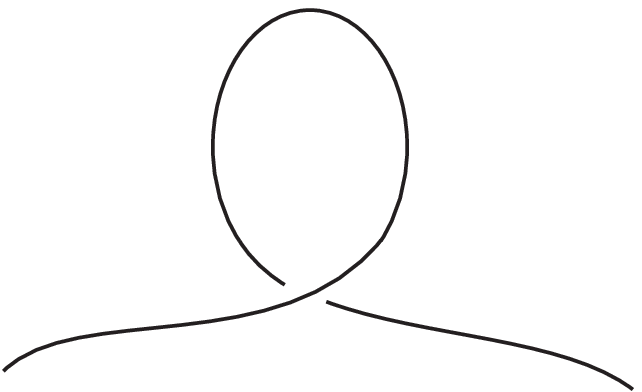}\end{overpic}
  \end{tabular}
  \caption{There can be no upper bound on $\delta([\gamma])$ for any knot type $[\gamma]$. Given any arc of a curve $\gamma$ of length $L$, as seen on the left, we can add a twist, as seen on the right. If $\gamma(s)$ lies directly over $\gamma(t)$ in this projection, then given any $\varepsilon > 0$ we can force $d(\gamma(s), \gamma(t); \R^3) < \varepsilon$ while $d(\gamma(s), \gamma(t); \gamma)$ remains constant.}
  \label{plie}
\end{figure}

In \cite{gromov} Gromov defines the distortion of a continuous curve $\gamma \colon [a,b] \longrightarrow \R^3$ as
\begin{equation*}
  \delta(\gamma)  := \sup_{s\ne t}\frac{d(\gamma(s), \gamma(t); \gamma)}{d(\gamma(s), \gamma(t); \R^3)}.
\end{equation*}
where $d(a, b; X)$ is the distance between $a$ and $b$ in the metric space $X$. We call the fraction inside the supremum the distortion quotient of the pair $(s,t)$, written $dq_{\gamma}(s,t)$. The distortion quotient becomes large when points that are some distance apart along $\gamma$ are close in space. For example, by adding a twist into a knot we can cause the distortion to grow arbitrarily large as seen in Figure~\vref{plie}. 
It follows that for any given knot type $[\gamma]$ and any given constant $M$, there exists a knot $\gamma \in [\gamma]$ so that $\delta(\gamma) > M$. Attempting to find the infimal distortion of a knot type is far more challenging. The distortion of a knot type, $\delta([\gamma]):=\inf_{\gamma \in [\gamma]} \delta(\gamma)$, is clearly bounded below by 1 for all knots, but finding the exact value of $\delta([\gamma])$ for any given knot type is very difficult. In fact, $\delta{[\gamma]}$ has only been computed when $[\gamma]$ is the unknot. As stated in \cite{gromov} and proved in \cite{KS} the distortion of the unknot is $\nicefrac{\pi}{2}$ and this value is achieved by the round circle. Gromov asks in \cite{gromov} whether or not there exists a universal upper bound on $\delta([\gamma])$. Though it has attracted considerable interest, Gromov's question has proved very difficult to answer.

It is natural to conjecture that no such bound exists. A proof would require two steps. We first need to find a topological invariant $\mathcal{X}$ that increases with knot complexity. We must then construct a lower bound on the distortion of a curve $\gamma$ in terms of $\mathcal{X}$. Using this lower bound, we can then exhibit a family of knots $\{\gamma_i\}$ so that the sequence $\{\delta([\gamma_i])\}$ diverges. Neither step seems easy to carry out.

Choosing a suitable invariant $\mathcal{X}$ requires care. It is known that a candidate topological invariant must increase without bound for some, but not all, families of knots: for a smooth knot $\gamma$ there exists an upper bound for the sequence $\{\delta([\gamma]), \delta([\gamma \# \gamma]), \delta([\gamma \#\gamma \# \gamma]), \ldots\}$\footnote{This was first pointed out in print by O'Hara \cite{ohara2}, who observed that the distortion of a number of tiny knots arranged around a large circle is independent of the number of knots.}. It follows that if the sequence $\{\mathcal{X}([\gamma]), \mathcal{X}([\gamma \# \gamma]), \mathcal{X}([\gamma \# \gamma \# \gamma]), \ldots\}$ diverges, then $\mathcal{X}$ will be of no use when defining a divergent sequence of lower bounds for distortion. This rules out many well known invariants such as crossing number, bridge number, and genus. 

The next task involves constructing a lower bound on the distortion of a curve in terms of the chosen topological invariant. This requires finding points $\gamma(s)$ and $\gamma(t)$ on the curve $\gamma$ for which $dq_{\gamma}(s,t)$ is large. Unfortunately, such points can be quite elusive. 

The main theorem of this paper, which summarizes the results of \cite{mullikin}, shows that length minimizing curves in a particular subset $U_C([\gamma])$ of $[\gamma]$ have the property that points with relatively large distortion quotient saturate arcs of positive total curvature.  This should make it easier to bound the distortion of curves in $U_C([\gamma])$ in terms of a topological invariant $\mathcal{X}([\gamma])$.

\begin{figure}[h]
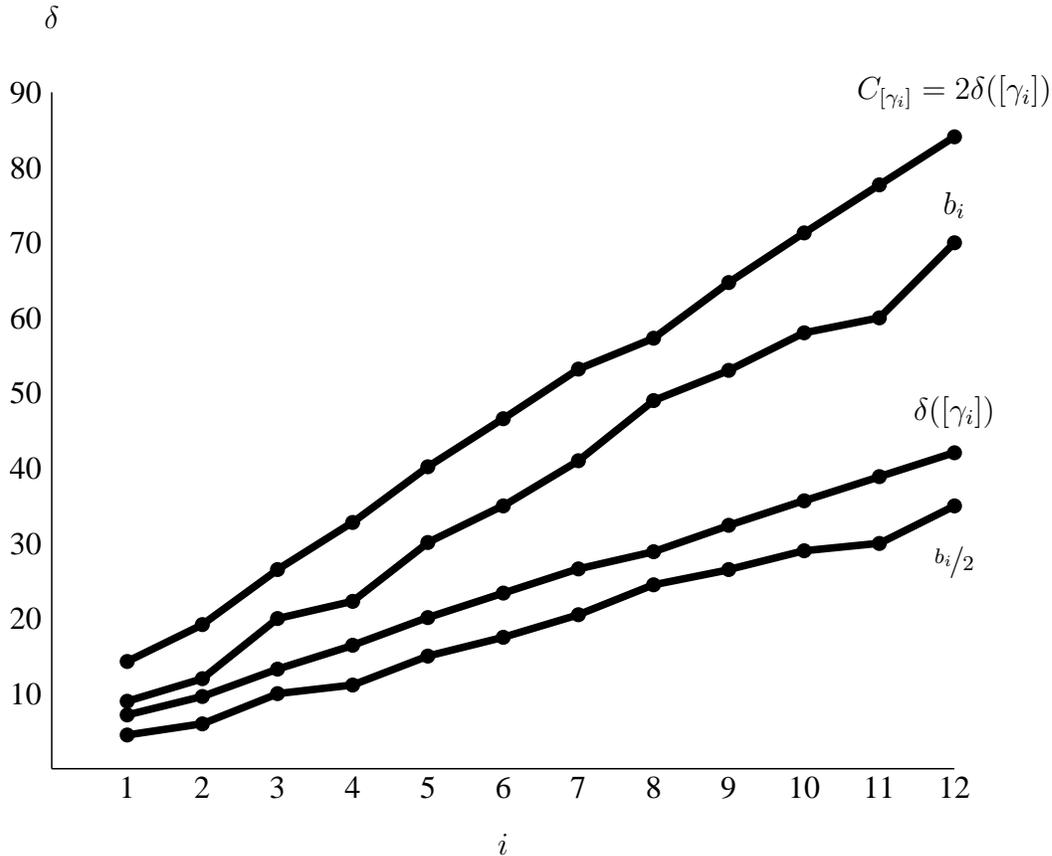

\begin{center}
          \begin{graph}(15,10)(0,-1)
            \roundnode{x}(0,9)[\graphnodesize{0.001}]
            \roundnode{y}(0,0)[\graphnodesize{0.001}]
            \roundnode{z}(12,0)[\graphnodesize{0.001}]
            \roundnode{xaxislabel}(6,-1)[\graphnodesize{0}]\autonodetext{xaxislabel}{$i$}
            \roundnode{yaxislabel}(0,10)[\graphnodesize{0}]\autonodetext{yaxislabel}{$\delta$}
            \roundnode{Clabel}(12,9)[\graphnodesize{0}]\autonodetext{Clabel}{$C_{[\gamma_i]} = 2\delta([\gamma_i])$}
            \roundnode{blabel}(12,7.5)[\graphnodesize{0}]\autonodetext{blabel}{$b_i$}        
            \roundnode{Blabel}(12,2.75)[\graphnodesize{0}]\autonodetext{Blabel}{$\nicefrac{b_i}{2}$}
            \roundnode{deltalabel}(12,4.75)[\graphnodesize{0}]\autonodetext{deltalabel}{$\delta([\gamma_i])$}
            \edge{x}{y}
            \edge{y}{z}
            \graphnodesize{0.2}
            %\roundnode{q1}(-0.35,1)[\graphnodesize{0}]\autonodetext{q1}{5}
            \roundnode{q2}(-0.35,1)[\graphnodesize{0}]\autonodetext{q2}{10}
            %\roundnode{q3}(-0.35,3)[\graphnodesize{0}]\autonodetext{q3}{15}
            \roundnode{q4}(-0.35,2)[\graphnodesize{0}]\autonodetext{q4}{20}
            %\roundnode{q5}(-0.35,5)[\graphnodesize{0}]\autonodetext{q5}{25}
            \roundnode{q6}(-0.35,3)[\graphnodesize{0}]\autonodetext{q6}{30}
            %\roundnode{q7}(-0.35,7)[\graphnodesize{0}]\autonodetext{q7}{35}
            \roundnode{q8}(-0.35,4)[\graphnodesize{0}]\autonodetext{q8}{40}
            %\roundnode{q9}(-0.35,9)[\graphnodesize{0}]\autonodetext{q9}{45}
            \roundnode{q9}(-0.35,5)[\graphnodesize{0}]\autonodetext{q9}{50}
            \roundnode{q10}(-0.35,6)[\graphnodesize{0}]\autonodetext{q10}{60}
            \roundnode{q11}(-0.35,7)[\graphnodesize{0}]\autonodetext{q11}{70}
            \roundnode{q12}(-0.35,8)[\graphnodesize{0}]\autonodetext{q12}{80}
            \roundnode{q13}(-0.35,9)[\graphnodesize{0}]\autonodetext{q13}{90}
            \roundnode{ly1}(0,0.5)[\graphnodesize{0}]
            \roundnode{ly2}(0,1)[\graphnodesize{0}]
            \roundnode{ly3}(0,1.5)[\graphnodesize{0}]
            \roundnode{ly4}(0,2)[\graphnodesize{0}]
            \roundnode{ly5}(0,2.5)[\graphnodesize{0}]
            \roundnode{ly6}(0,3)[\graphnodesize{0}]
            \roundnode{ly7}(0,3.5)[\graphnodesize{0}]
            \roundnode{ly8}(0,4)[\graphnodesize{0}]
            \roundnode{ly9}(0,4.5)[\graphnodesize{0}]
            \roundnode{ly10}(0,5)[\graphnodesize{0}]
            \roundnode{ly11}(0,5.5)[\graphnodesize{0}]
            \roundnode{ly12}(0,6)[\graphnodesize{0}]
            \roundnode{ly13}(0,6.5)[\graphnodesize{0}]
            \roundnode{ly14}(0,7)[\graphnodesize{0}]
            \roundnode{ly15}(0,7.5)[\graphnodesize{0}]
            \roundnode{ly16}(0,8)[\graphnodesize{0}]
            \roundnode{ly17}(0,8.5)[\graphnodesize{0}]
            \roundnode{ly18}(0,9)[\graphnodesize{0}]
            \roundnode{ry1}(12,0.5)[\graphnodesize{0}]
            \roundnode{ry2}(12,1)[\graphnodesize{0}]
            \roundnode{ry3}(12,1.5)[\graphnodesize{0}]
            \roundnode{ry4}(12,2)[\graphnodesize{0}]
            \roundnode{ry5}(12,2.5)[\graphnodesize{0}]
            \roundnode{ry6}(12,3)[\graphnodesize{0}]
            \roundnode{ry7}(12,3.5)[\graphnodesize{0}]
            \roundnode{ry8}(12,4)[\graphnodesize{0}]
            \roundnode{ry9}(12,4.5)[\graphnodesize{0}]
            \roundnode{ry10}(12,5)[\graphnodesize{0}]
            \roundnode{ry11}(12,5.5)[\graphnodesize{0}]
            \roundnode{ry12}(12,6)[\graphnodesize{0}]
            \roundnode{ry13}(12,6.5)[\graphnodesize{0}]
            \roundnode{ry14}(12,7)[\graphnodesize{0}]
            \roundnode{ry15}(12,7.5)[\graphnodesize{0}]
            \roundnode{ry16}(12,8)[\graphnodesize{0}]
            \roundnode{ry17}(12,8.5)[\graphnodesize{0}]
            \roundnode{ry18}(12,9)[\graphnodesize{0}]
            \roundnode{p1}(1,0.7158382)
            \roundnode{p2}(2,0.9624142)
            \roundnode{p3}(3,1.3271537) 
            \roundnode{p4}(4,1.6421985)
            \roundnode{p5}(5,2.0126675)
            \roundnode{p6}(6,2.3375325)
            \roundnode{p7}(7,2.6617379)
            \roundnode{p8}(8,2.8904931)
            \roundnode{p9}(9, 3.2382402)
            \roundnode{p10}(10,3.5667338)
            \roundnode{p11}(11,3.8863978)
            \roundnode{p12}(12,4.2067961)
            \roundnode{B1}(1,.9)
            \roundnode{B2}(2,1.2)
            \roundnode{B3}(3,2) 
            \roundnode{B4}(4,2.23)
            \roundnode{B5}(5,3.01)
            \roundnode{B6}(6,3.5)
            \roundnode{B7}(7,4.1)
            \roundnode{B8}(8,4.9)
            \roundnode{B9}(9, 5.3)
            \roundnode{B10}(10,5.8)
            \roundnode{B11}(11,6)
            \roundnode{B12}(12,7)
            \roundnode{b1}(1,0.45)
            \roundnode{b2}(2,0.6)
            \roundnode{b3}(3,1) 
            \roundnode{b4}(4,1.115)
            \roundnode{b5}(5,1.5005)
            \roundnode{b6}(6,1.75)
            \roundnode{b7}(7,2.05)
            \roundnode{b8}(8,2.45)
            \roundnode{b9}(9, 2.65)
            \roundnode{b10}(10,2.9)
            \roundnode{b11}(11,3)
            \roundnode{b12}(12,3.5)
            \roundnode{c1}(1,1.43)
            \roundnode{c2}(2,1.92)
            \roundnode{c3}(3,2.65) 
            \roundnode{c4}(4,3.28)
            \roundnode{c5}(5,4.02)
            \roundnode{c6}(6,4.66)
            \roundnode{c7}(7,5.32)
            \roundnode{c8}(8,5.73)
            \roundnode{c9}(9,6.47)
            \roundnode{c10}(10,7.13)
            \roundnode{c11}(11,7.77)
            \roundnode{c12}(12,8.41)
            \roundnode{r1}(1,-0.25)[\graphnodesize{0}]\autonodetext{r1}{1}
            \roundnode{r2}(2,-0.25)[\graphnodesize{0}]\autonodetext{r2}{2}
            \roundnode{r3}(3,-0.25)[\graphnodesize{0}]\autonodetext{r3}{3}
            \roundnode{r4}(4,-0.25)[\graphnodesize{0}]\autonodetext{r4}{4}
            \roundnode{r5}(5,-0.25)[\graphnodesize{0}]\autonodetext{r5}{5}
            \roundnode{r6}(6,-0.25)[\graphnodesize{0}]\autonodetext{r6}{6}
            \roundnode{r7}(7,-0.25)[\graphnodesize{0}]\autonodetext{r7}{7}
            \roundnode{r8}(8,-0.25)[\graphnodesize{0}]\autonodetext{r8}{8}
            \roundnode{r9}(9,-0.25)[\graphnodesize{0}]\autonodetext{r9}{9}
            \roundnode{r10}(10,-0.25)[\graphnodesize{0}]\autonodetext{r10}{10}
            \roundnode{r11}(11,-0.25)[\graphnodesize{0}]\autonodetext{r11}{11}
            \roundnode{r12}(12,-0.25)[\graphnodesize{0}]\autonodetext{r12}{12}
            \roundnode{bx1}(1,0)[\graphnodesize{0}]
            \roundnode{bx2}(2,0)[\graphnodesize{0}]
            \roundnode{bx3}(3,0)[\graphnodesize{0}]
            \roundnode{bx4}(4,0)[\graphnodesize{0}]
            \roundnode{bx5}(5,0)[\graphnodesize{0}]
            \roundnode{bx6}(6,0)[\graphnodesize{0}]
            \roundnode{bx7}(7,0)[\graphnodesize{0}]
            \roundnode{bx8}(8,0)[\graphnodesize{0}]
            \roundnode{bx9}(9,0)[\graphnodesize{0}]
            \roundnode{bx10}(10,0)[\graphnodesize{0}]
            \roundnode{bx11}(11,0)[\graphnodesize{0}]
            \roundnode{bx12}(12,0)[\graphnodesize{0}]
            \roundnode{tx1}(1,9)[\graphnodesize{0}]
            \roundnode{tx2}(2,9)[\graphnodesize{0}]
            \roundnode{tx3}(3,9)[\graphnodesize{0}]
            \roundnode{tx4}(4,9)[\graphnodesize{0}]
            \roundnode{tx5}(5,9)[\graphnodesize{0}]
            \roundnode{tx6}(6,9)[\graphnodesize{0}]
            \roundnode{tx7}(7,9)[\graphnodesize{0}]
            \roundnode{tx8}(8,9)[\graphnodesize{0}]
            \roundnode{tx9}(9,9)[\graphnodesize{0}]
            \roundnode{tx10}(10,9)[\graphnodesize{0}]
            \roundnode{tx11}(11,9)[\graphnodesize{0}]
            \roundnode{tx12}(12,9)[\graphnodesize{0}]
            \edge{p1}{p2}[\graphlinewidth{.1}]
            \edge{p2}{p3}[\graphlinewidth{.1}]
            \edge{p3}{p4}[\graphlinewidth{.1}]
            \edge{p4}{p5}[\graphlinewidth{.1}]
            \edge{p5}{p6}[\graphlinewidth{.1}]
            \edge{p6}{p7}[\graphlinewidth{.1}]
            \edge{p7}{p8}[\graphlinewidth{.1}]
            \edge{p8}{p9}[\graphlinewidth{.1}]
            \edge{p9}{p10}[\graphlinewidth{.1}]
            \edge{p10}{p11}[\graphlinewidth{.1}]
            \edge{p11}{p12}[\graphlinewidth{.1}]
            \edge{c1}{c2}[\graphlinewidth{.1}]
            \edge{c2}{c3}[\graphlinewidth{.1}]
            \edge{c3}{c4}[\graphlinewidth{.1}]
            \edge{c4}{c5}[\graphlinewidth{.1}]
            \edge{c5}{c6}[\graphlinewidth{.1}]
            \edge{c6}{c7}[\graphlinewidth{.1}]
            \edge{c7}{c8}[\graphlinewidth{.1}]
            \edge{c8}{c9}[\graphlinewidth{.1}]
            \edge{c9}{c10}[\graphlinewidth{.1}]
            \edge{c10}{c11}[\graphlinewidth{.1}]
            \edge{c11}{c12}[\graphlinewidth{.1}]
            \edge{b1}{b2}[\graphlinewidth{.1}]
            \edge{b2}{b3}[\graphlinewidth{.1}]
            \edge{b3}{b4}[\graphlinewidth{.1}]
            \edge{b4}{b5}[\graphlinewidth{.1}]
            \edge{b5}{b6}[\graphlinewidth{.1}]
            \edge{b6}{b7}[\graphlinewidth{.1}]
            \edge{b7}{b8}[\graphlinewidth{.1}]
            \edge{b8}{b9}[\graphlinewidth{.1}]
            \edge{b9}{b10}[\graphlinewidth{.1}]
            \edge{b10}{b11}[\graphlinewidth{.1}]
            \edge{b11}{b12}[\graphlinewidth{.1}]  
            \edge{B1}{B2}[\graphlinewidth{.1}]
            \edge{B2}{B3}[\graphlinewidth{.1}]
            \edge{B3}{B4}[\graphlinewidth{.1}]
            \edge{B4}{B5}[\graphlinewidth{.1}]
            \edge{B5}{B6}[\graphlinewidth{.1}]
            \edge{B6}{B7}[\graphlinewidth{.1}]
            \edge{B7}{B8}[\graphlinewidth{.1}]
            \edge{B8}{B9}[\graphlinewidth{.1}]
            \edge{B9}{B10}[\graphlinewidth{.1}]
            \edge{B10}{B11}[\graphlinewidth{.1}]
            \edge{B11}{B12}[\graphlinewidth{.1}]             
          \end{graph}
        \end{center}
        \caption{Here we see a graph representing a divergent sequence $\{b_i\}$ of lower bounds on the distortion of certain curves in $U_{C_{\gamma_i}}([\gamma_i])$. If we define $C_{\gamma_i} = 2\delta([\gamma_i])$ then this will also yield a divergent sequence $\{\nicefrac{b_i}{2}\}$ of lower bounds on $\delta([\gamma_i])$.}
        \label{HolySmokes}
\end{figure}

The set $U_C([\gamma])$ has several membership requirements in addition to the fact that $U_C([\gamma]) \subset [\gamma]$. For technical reasons we assume all curves in $U_C([\gamma])$ have finite total curvature. Second, we require that a constant $C$ bounds the distortion of knots in $U_C([\gamma])$ from above. If we choose $C$ to be a constant multiple of $\delta([\gamma])$, this bound will allow us to relate lower bounds for distortion of curves in $U_C([\gamma])$ to lower bounds on $\delta([\gamma])$. Hence any divergent sequence of lower bounds $\{b_i\}$ on the distortion of knots in, say, $\{U_{2\delta([\gamma])}([\gamma_i])\}$ will provide a negative answer to Gromov's question. This is illustrated in Figure~\vref{HolySmokes}.

\begin{figure}[htbp]
  \centerline{\begin{overpic}{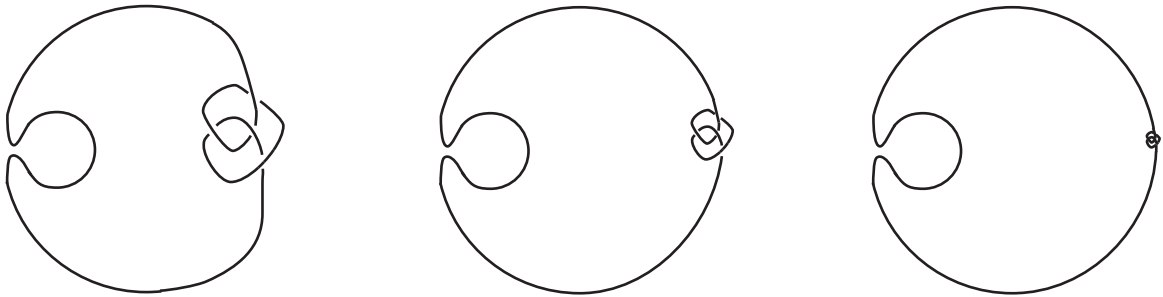}\end{overpic}}
  \caption{Here are three representatives from a sequence of figure eight knots with decreasing length that converges to the unknot. The pinch on the left illustrates the concern that every member of the sequence of curves may have maximum distortion realized by a pair of points a constant distance apart even though the knotted region is shrinking to a point.}
  \label{limit_danger}
\end{figure}

The last requirement for membership in $U_C([\gamma])$ stems from our interest in studying curves of minimum length. Since distortion is scale invariant the set of length minimizers in $U_C([\gamma])$ will be empty unless we add another constraint. Hence we must fix a scale for curves in $U_C([\gamma])$ to prevent any sequence of curves $\{\gamma_k\} \subset U_C([\gamma])$ approaching a curve $\gamma_0$ of infimal length from shrinking to a point. However, we still need to be concerned that $\gamma_0$ might not be in $[\gamma]$. After all, knotted regions can pull tight when decreasing length as seen in Figure~\vref{limit_danger}.  To prevent this, we fix scale in a carefully chosen way. In \cite{KS} Kusner and Sullivan define the $b$-distortion thickness of a curve $\gamma$ as
\begin{equation*}
  \tau_b(\gamma) = \inf_{dq_{\gamma}(p,q) \geq b} \|\gamma(p) - \gamma(q)\|,
\end{equation*}
the infimal distance between pairs of points on $\gamma$ with distortion quotient at least $b$. We require that $\tau_{\delta([\gamma])}(\gamma) \geq 1$ for curves $\gamma \in U_C([\gamma])$.

We are now prepared to state the main theorem. For now we think of a ``$\delta([\gamma])$-drc" as a pair of points $(s,t)$ for which $dq_{\gamma}(s,t) = \delta([\gamma])$. We will define this term precisely in Definition~\vref{kdrc}.

\setcounter{thm}{0}
  \begin{thm} [Main Theorem] Let $U_C([\gamma])$ be the set of all finite total curvature curves $\gamma$ in $[\gamma]$, with distortion $\delta(\gamma) < C$ and distortion thickness $\tau_{\delta([\gamma])} \geq 1$.

 Then any open interval on a curve of minimum length in $U_C([\gamma])$ is either a straight line segment or contains an endpoint of a $\delta([\gamma])$-drc.
\end{thm}
\setcounter{thm}{0}

The proof of the main theorem will be a proof by contradiction. Suppose that we have a curve $\gamma\colon[a,b] \longrightarrow \R^3$ of minimum length in $U_C([\gamma])$ and an arc $\gamma((c,d))$ on $\gamma$ with positive total curvature that does not contain an endpoint of a $\delta([\gamma])$-drc. We show that it is possible to decrease the length of the arc $\gamma((c,d))$ to obtain a new curve that is also a member of $U_C([\gamma])$. We have then reached a contradiction since $\gamma$ is a curve of minimum length in $U_C([\gamma])$. The proof requires two propositions.

\begin{itemize}
    \item[(1)] There exists $\varepsilon > 0$ so that $dq(s,t) < \delta([\gamma]) - \varepsilon$ for all points $(s,t) \in (c,d) \times [a,b] \cup [a,b] \times (c,d)$. (Proposition~\vref{prop2}) 

  \item[(2)] The arc $\gamma((c,d))$ can be changed in a length decreasing way so that the increase of the distortion quotient on $(c,d) \times [a,b] \cup [a,b] \times (c,d)$ is less than $\varepsilon$ and the change of the distortion quotient on the remainder of $[a,b] \times [a,b]$ is nonpositive. Hence the distortion of $\gamma$ is not increased.  (Proposition~\vref{prop1})
\end{itemize}

We point out that if $dq_{\gamma}$ could be extended to a continuous function on $[a,b] \times [a,b]$ then (1) would be immediate and (2) would be much easier to prove. However, this is not always possible for finite total curvature curves as we will see in section 3.

\section{Definitions and Background}

The proof of our main theorem makes use of several standard results from analysis as well as few results more directly related to distortion. We will be dealing with the set of curves with finite total curvature in the sense of Milnor in \cite{milnor}. Given a polygon $p$ with vertices $v_0, \ldots, v_n$, we can define a sequence of angles $\alpha_0, \ldots, \alpha_n$ where $\alpha_i$ is the exterior angle between the two edges incident to $v_i$. In \cite{milnor}, Milnor defines the total curvature of the polygon $p$ to be $\kappa_M(p) := \sum_{i = 0}^n \alpha_i$. He then defines the total curvature of a curve $\gamma$ to be $\sup_{p \in \text{Pol}(\gamma)}\kappa_M(p)$, where $\text{Pol}(\gamma)$ denotes the set of all polygons inscribed in $\gamma$. It is also shown in \cite{milnor} that if $\gamma \in C^2$, then $\kappa_M(\gamma)$ agrees with the standard definition for total curvature. There are many interesting properties of the class of finite total curvature curves, denoted $\FTC$, which are explored in detail in both \cite{Sullivan-FTC} and \cite{mullikin}. In particular, if $\gamma \in \FTC$, then $\gamma$ has one-sided derivatives everywhere. Therefore we can define both the left and right tangent indicatrices $T_{\gamma}^{\pm}\colon [a,b] \longrightarrow S^1$by
\begin{equation*}
  T_{\gamma}^{\pm}(s) := \lim_{x \rightarrow s^{\pm}}\frac{\gamma(x) - \gamma(s)}{\|\gamma(x) - \gamma(s)\|}.
\end{equation*}
It can be shown that each of these curves is of bounded variation and therefore, by the structure theorem for BV functions, they possess a derivative in the sense of the norm of their derivative is a Radon measure $\mathcal{K}$. We find an explicit construction of this measure $\mathcal{K}$ in \cite{mullikin}. The measure theoretic properties of $\mathcal{K}$ which are required for the proof are given in Lemma~\vref{measuregoodness} whose proof we move to the appendix.

The distortion of a curve increases as points on the curve become closer together in space while remaining relatively far apart in arclength. For arcs of curves that are nearly straight the value of the distortion quotient for pairs on this arc must remain fairly close to 1. Indeed Denne and Sullivan show in \cite{DS} the following relationship between total curvature and distortion.

\begin{lemma}(Denne/Sullivan) \label{DSbound} Let $\gamma \colon [a,b] \longrightarrow \R^3$ be any finite total curvature curve and let $\gamma([c,d])$ be any arc of $\gamma$ with total curvature $\kappa \leq \pi$. Then $dq_{\gamma}(s,t) \leq \sec{\nicefrac{\kappa}{2}}$ for all $(s,t) \in [c,d] \times [c,d]$ .
\end{lemma}

We omit the proof of Lemma~\ref{DSbound} which appears in \cite{DS} and \cite{mullikin}.

\section{On the Discontinuity of $dq_{\gamma} \text{ for } \gamma \in$ FTC}

Let $\gamma \colon [a,b] \longrightarrow \R^3$ be a continuous embedded curve. Then since $\gamma$ is continuous, the functions $d(\cdot, \cdot; \gamma)$ and $d(\cdot, \cdot; \R^3)$, restricted to points on $\gamma$, are each continuous functions. So their ratio, the distortion quotient, is continuous whenever the denominator is nonzero. Indeed, $dq_{\gamma}$ is only defined for points off the diagonal of $[a,b] \times [a,b]$.

Regrettably, if $\gamma \in \FTC$ then it may not be possible to define $dq_{\gamma}$ on the diagonal in a way that will result in a continuous function. While curves in $\FTC$ have one-sided tangents everywhere, the right and left tangent may not be equal at a given point $s_0$. This can cause a non-removable discontinuity of $dq_{\gamma}$ at $(s_0, s_0)$. The following example illustrates such a scenario.

\begin{figure}[h]
  \centerline{\begin{overpic}{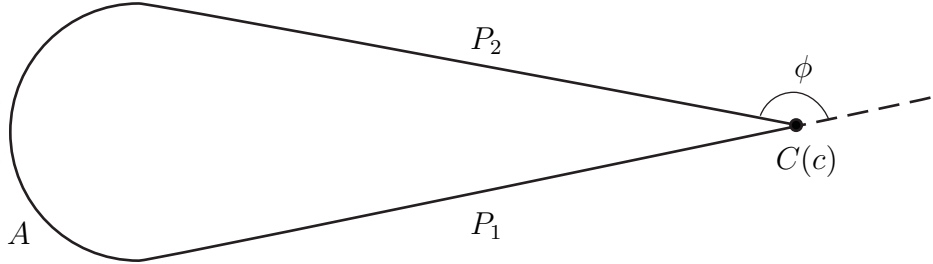}
  \put(0,2){$A$}
  \put(50,3){$P_1$}
  \put(50,23){$P_2$}
  \put(85,20){$\phi$}
  \put(83,10){$C(c)$}
  \end{overpic}}
  \caption{This example illustrates a discontinuity of the distortion quotient. The discontinuity of the tangent curve at the point corresponding to the corner point causes a jump in the distortion quotient.} 
    \label{comet}
  \end{figure} 

 Let $\mathcal{C} \colon [a,b] \longrightarrow \R^2$ denote the comet shaped curve in Figure \ref{comet} consisting of an arc $A$ of a circle and a polygonal section, consisting of two line segments $P_1$ and $P_2$, with exterior angle $\phi$. Let $c \in [a,b]$ so that $\mathcal{C}(c)$ is the corner point. A calculation shows that $\delta(\mathcal{C}) = \sec{(\phi/2)}$ and $dq_{\mathcal{C}}(s,t) = \delta(\mathcal{C})$ on the continuum of points 
 \begin{equation*}
 E = \{(s,t) \in [a,b] \times [a,b] : \mathcal{C}(s) \in P_1, \text{ } \mathcal{C}(t) \in P_2\text{, and }d(\mathcal{C}(s), \mathcal{C}(c), \mathcal{C}) = d(\mathcal{C}(t), \mathcal{C}(c), \mathcal{C})\}.
 \end{equation*}

\begin{figure}[htbp]
  \centerline{\begin{overpic}{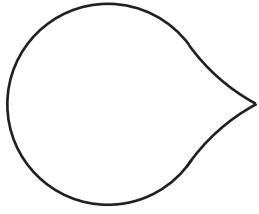}\end{overpic}}
  \caption{The distortion of this curve is realized as a limit of distortion quotients of pairs of symmetric points converging to the corner point. The parameter space $[a,b] \times [a,b]$ has no point $(s,t)$ that realizes the distortion.}
   \label{dragon_tooth_pic}
\end{figure}

Now, if $\{(s_i, t_i)\} \subset E$ is a sequence so that $(s_i, t_i) \rightarrow (c,c)$ as $i \rightarrow \infty$, then 
\begin{equation*}
  \lim_{i \rightarrow \infty}dq_{\mathcal{C}}(s_i,t_i) = \sec \frac{\phi}{2}.
\end{equation*}
On the other hand, for all points $s$ so that $\mathcal{C}(s)$ is a point on either one of theline segments, we have $dq_{\mathcal{C}}(s,c)~=~1$, hence
\begin{equation*}
  \lim_{s \rightarrow c}dq_{\mathcal{C}}(s,c) = 1.
\end{equation*}
We can readily see that the discontinuity along the diagonal of $[a,b] \times [a,b]$ is not removable. In fact, since the distortion is defined as a supremum, it may be the case that for a general curve $\gamma \in \FTC$, there exists a sequence of points $(s_i, t_i)$ converging to a point $(s,s)$ so that $dq_{\gamma}(s_i, t_i)~\rightarrow~\delta(\gamma)$, \hspace{0.1in}$dq_{\gamma}(s,t)~<~\delta(\gamma)$ everywhere $dq_{\gamma}$ is defined, and $dq_{\gamma}$ has a non-removable discontinuity at $(s,s)$. Curves of this type have no distortion realizing chord. An example of such a curve is the \emph{Dragon's tooth curve} in Figure~\ref{dragon_tooth_pic}.

This curve, $\mathcal{D} \colon [a,b] \longrightarrow \R^2$, is formed by three circle arcs: one of small radius $r$ and two of larger radius $R$. The large circle arcs are connected to the smaller circle arc so that their tangents agree at the points of intersection. The large circles meet at a corner with exterior angle $\phi~<~\pi$. The distortion of a circle is $\pi/2$, so if $\mathcal{D}(s)$ and $\mathcal{D}(t)$ are points on the same circle arc then $dq_{\mathcal{D}}(s,t)~\leq~\pi/2$. If $\mathcal{D}(s)$ lies on one of the large circle arcs and if $\mathcal{D}(t)$ lies on the small circle arc, then $dq_{\mathcal{D}}(s,t) \rightarrow 1$ as $d(\mathcal{D}(s), \mathcal{D}(t); \R^2) \rightarrow 0$. Therefore, the distortion quotient cannot grow large because the denominator becomes small. On the other hand, the numerator is no larger than the length of the two arcs. So $dq_{\mathcal{D}}(s,t)$ is uniformly bounded above for all points $(s,t)$ for which $\mathcal{D}(s)$ and $\mathcal{D}(t)$ lie on circle arcs of different radius. 

Let $\mathcal{D}(c)$ denote the corner point. Let $\{(s_i, t_i)\}$ be a sequence so that $\mathcal{D}(s_i)$ is on one large circle arc, $\mathcal{D}(t_i)$ is on the other large circle arc and $d(\mathcal{D}(s_i), \mathcal{D}(c); \mathcal{D}) = d(\mathcal{D}(t_i), \mathcal{D}(c); \mathcal{D})$ for all $i$. If $(s_i, t_i)\rightarrow (c,c)$ then
\begin{equation*}
  \lim_{i \rightarrow \infty}dq_{\mathcal{D}}(s_i, t_i) = \sec\frac{\phi}{2}.
\end{equation*}
Forcing $\phi$ to be close to $\pi$ in the construction of $\mathcal{D}$ forces $dq_{\mathcal{D}}(s,t) < \sec\phi/2$ for all $(s,t)$ in the domain of $dq_{\mathcal{D}}$. So, the distortion of $\mathcal{D}$ is evaluated as a limit of $dq_{\mathcal{D}}(s,t)$ as $(s,t)$ converges to a point on the diagonal that is a non-removable discontinuity. 

To deal with examples like this, the definition of a distortion realizing chord will need to be extended. Just considering pairs of points $(s,t)$ so that $dq_{\gamma}(s,t) = \delta(\gamma)$ is not sufficient. After all, there may be no such pairs. We start with a definition.

\begin{Def}
For any  $\gamma \colon [a,b] \longrightarrow \R^3$ in $\FTC$, we define the function $D_{\gamma} \colon [a,b] \longrightarrow \R$ by 
\begin{equation*}
  D_{\gamma}(s) \colon =\sup_t dq_{\gamma}(s,t)
\end{equation*}
\end{Def}

\begin{lemma}
  If $\gamma \colon [a,b] \longrightarrow \R^3$ is a curve with finite total curvature parametrized by arclength, then
\begin{equation*}
  D_{\gamma}(s) = \sup_{t} \frac{d(\gamma(s),\gamma(t); \gamma)}{d(\gamma(s),\gamma(t);\R^3)}  =\max_{t} \frac{d(\gamma(s),\gamma(t); \gamma)}{d(\gamma(s),\gamma(t);\R^3)}.
  \end{equation*}
  \label{D=max}
\end{lemma}

\begin{proof} 
It can be shown that one-sided derivatives exist everywhere for curves in $\FTC$ (see \cite{mullikin}), hence
\begin{equation*}
  \lim_{t\rightarrow s^+}\frac{\gamma(t) - \gamma(s)}{t - s}
\end{equation*}
exists. We know, by the continuity of the norm function $\| \cdot \|$,  that
\begin{equation*}
    \lim_{t\rightarrow s^+}\left\|\frac{\gamma(t) - \gamma(s)}{t - s}\right\| = \lim_{t\rightarrow s^+}\frac{\|\gamma(t) - \gamma(s)\|}{|t - s|}
\end{equation*}
exists and equals 1 since the curve is arclength parametrized. So,
\begin{equation*}
  \lim_{t\rightarrow s^+}\frac{|t - s|}{\|\gamma(t) - \gamma(s)\|}=1.
\end{equation*}

But this is exactly the one-sided limit of the distortion quotient function $dq_{\gamma}(s,t)$ if we leave $s$~fixed. A similar computation shows that the limit from the left exists as well and is also equal to~1. Therefore, since the function $dq_{\gamma}(s,t)$ with $s$ fixed is a continuous function on the compact set $[a,b]$, it follows that it achieves its maximum at a point in $[a,b]$.
\end{proof}

In the case of the Dragon's tooth $\mathcal{D}$, the function $D_{\mathcal{D}}(s)$ is defined at the corner point, but it is less than the distortion of the curve. However, we do have $D_{\mathcal{D}}(s_i) \rightarrow \delta(\mathcal{D})$ as $s_i \rightarrow s$. We will replace $D_{\gamma}$ so that our new function will achieve its supremal value on $\gamma$.

\begin{Def} If $\gamma \colon [a,b] \longrightarrow \R^3$ is a curve with finite total curvature parametrized by arclength, then we define the \textbf{distortion shadow}, denoted $\overline{D}_{\gamma}(s)$,  to be the upper envelope of $D_{\gamma}(s)$. That is
\begin{equation*}
  \overline{D}_{\gamma}(s) \colon = \inf_{\varepsilon > 0} \sup_{|s-x|<\varepsilon} \left(\max_{t} \frac{d(\gamma(x),\gamma(t); \gamma)}{d(\gamma(x), \gamma(t); \R^3)}\right).
\end{equation*}
\end{Def}

We are now in a position to make the definition of a distortion realizing chord precise. The Dragon's tooth example shows that we will need to define the distortion quotient along the diagonal in a way that guarantees that the value of $dq_{\gamma}(s,s)$ is defined to be the largest limiting value of all sequences approaching it. For this we again use the upper envelope.  

\begin{Def}\label{kdrc}
  Let $\gamma \colon [a,b] \longrightarrow \R^3$ be a finite total curvature curve parametrized by arclength. If $(s,t)$ is such that $s \ne t$ and
\begin{equation*}
  \overline{dq_{\gamma}}(s,t)  := \inf_{\varepsilon >0} \sup_{\|(x,y) - (s,t)\|< \varepsilon} \frac{d(\gamma(x),\gamma(y); \gamma)}{d(\gamma(x),\gamma(y); \R^3)} \geq k,
\end{equation*}
then we say $(s,t)$ (or the chord with endpoints $\gamma(s)$ and $\gamma(t)$) is a $k$\textbf{-distortion realizing chord} ($k$\textbf{-drc}).
\end{Def}

\begin{lemma}
  \label{one}
  Let $\gamma \colon [a,b] \longrightarrow \R^3$ be a finite total curvature curve parametrized by arclength. Then $\overline{D}_{\gamma}(s) \geq k$ if and only if there exists a value $t$ so that the chord with endpoints $\gamma(s)$ and $\gamma(t)$ is a $k$-drc.
\end{lemma}

\begin{proof}  
First assume that there is a value $t$ so that the point $(s,t)$ defines a $k$-drc. We know that there exists a sequence $\{(s_i, t_i)\}$ so that $(s_i, t_i) \rightarrow (s,t)$ and $dq_{\gamma}(s_i, t_i)\rightarrow K \geq k$. Since the function $D_{\gamma}(s_i)$ computes the maximum over all values of $t$, it is evident that $D_{\gamma}(s_i) \geq dq_{\gamma}(s_i, t_i)$. Furthermore, since the upper envelope is upper-semicontinuous, $\overline{D}_{\gamma}(s_i) \geq D_{\gamma}(s_i)$. So our string of inequalities then becomes
\begin{equation*}
    \overline{D}_{\gamma}(s) \geq \lim \overline{D}_{\gamma}(s_i) \geq \lim D_{\gamma}(s_i) \geq \lim dq_{\gamma}(s_i, t_i) = K \geq k.
\end{equation*}

The reverse implication is immediate from the definitions. Let $\{s_i\} \subset [a,b]$ be a sequence so that $D_{\gamma}(s_i) \rightarrow K \geq k$. For each $s_i$ there exists a $t_i$ so that $D(s_i) = dq_{\gamma}(s_i, t_i)$ by Lemma~\ref{D=max}. Furthermore, since the points $(s_i, t_i)$ are elements of the compact set $[a,b]\times[a,b]$ we may assume, by restricting to a subsequence, that $(s_i, t_i) \rightarrow (s,t)$ for some $(s,t) \in [a,b]\times[a,b]$. It remains to show $\overline{dq}_{\gamma}(s,t) \geq k$. Indeed, for any $\varepsilon > 0$
\begin{equation*}
  \sup_{\|(x,y) - (s,t)\|<\varepsilon}dq(x,y)\geq k
\end{equation*}
since the set on which the supremum is taken contains infinitely many elements of the sequence $\{(s_i, t_i)\}$. Therefore,
\begin{equation*}
\overline{dq}_{\gamma}(s,t) =   \inf_{\varepsilon > 0} \sup_{\|(x,y) - (s,t)\|<\varepsilon}dq(s,t) \geq k.
\end{equation*}
So the pair $(s,t)$ defines a $k$-drc.
\end{proof}

By defining $k$-drc's in terms of an upper-semicontinuous function we have made the proof of the first step towards the main theorem relatively simple. Indeed, upper-semicontinuous functions achieve their maximum value on compact sets. Therefore, if a closed arc contains no endpoints of $k$-drcs, then the maximum value of $\overline{dq_{\gamma}}$ is bounded away from $k$ by a positive quantity $\varepsilon$. Hence, the distortion quotient is bounded away from $k$ by a positive quantity. We make this rigorous below.

\begin{proposition} \label{prop2}Let $\gamma \colon [a,b] \longrightarrow \R^3$ be a curve of finite total curvature. If $(c,d) \subset [a,b]$ so that $\overline{dq_{\gamma}}(s,t) < k$ for all $(s,t) \in (c,d) \times [a,b] \cup [a,b] \times (c,d)$ then, for any subinterval $[p,q]~\subset~(c,d)$ there exists an $\varepsilon > 0$ so that $dq_{\gamma}(s,t) \leq k - \varepsilon$ for all $(s,t) \in (p,q) \times [a,b] \cup [a,b] \times (p,q)$
\end{proposition}

\begin{proof} 
By assumption the interval $(c,d)$ is free from endpoints of $k$-drc's. Then using Lemma~\vref{one}, it follows that for every $s \in (c,d)$ the value of $\overline{D}_{\gamma}(s) < k$. Let $[p,q]$ be any closed interval subset of $(c,d)$. Since $\overline{D}_{\gamma}(s)$ is upper-semicontinuous and bounded, $\overline{D}_{\gamma}(s)$ has a maximum $M$ on $[p,q]$. We can let $\varepsilon = k - M$.
\end{proof}

\section{Decreasing Length Without Increasing Distortion}

Now that we know that arcs of $\gamma$ which do not contain $k$-drc's have distortion bounded away from $\delta(\gamma)$ by some positive quantity $\varepsilon$, it remains to show the length of an arc can be decreased while changing the distortion quotient by an amount smaller than $\varepsilon$. This will complete the second step in the proof of the main theorem outlined in Section 2.1.

\begin{proposition}\label{prop1} Let $\gamma \colon [a,b] \longrightarrow \R^3$ be a finite total curvature curve and suppose $\gamma((c,d))$ is any arc with nonzero total curvature. Then, given $\varepsilon > 0$, we can replace $\gamma((c,d))$ with an arc $P((c,)])$ of shorter length so that $dq_P(s,t) - dq_{\gamma}(s,t) < \varepsilon$ for all $(s,t)~\in~(c,d)~\times~[a,b] \cup [a,b] \times (c,d)$.
\end{proposition}

\begin{proof} Let $\gamma \colon [a,b] \longrightarrow \R^3$ be a curve parametrized by arclength and suppose that $\gamma((c,d))$ is an arc of $\gamma$ of length $L = d - c$ with nonzero total curvature. Let $\varepsilon > 0$ be given. This proof has several cases. We will first assume that we can find a subarc of $\gamma((c,d))$ whose total curvature is small. Without loss of generality we assume this is the arc $\gamma((c,d))$. We will then modify only a subarc $\gamma([\hat{c}, \hat{d}])$ of $\gamma((c,d))$. To prove that distortion has not greatly increased, we will then examine the change in the distortion quotient. Pairs of points $s$ and $t$ which are relatively close will require a different argument than pairs where $s$ and $t$ are relatively far apart. Figure~\ref{reference} may help the reader to understand the different components of the proof. 

\begin{figure}[ht!]
    \centerline{\begin{overpic}{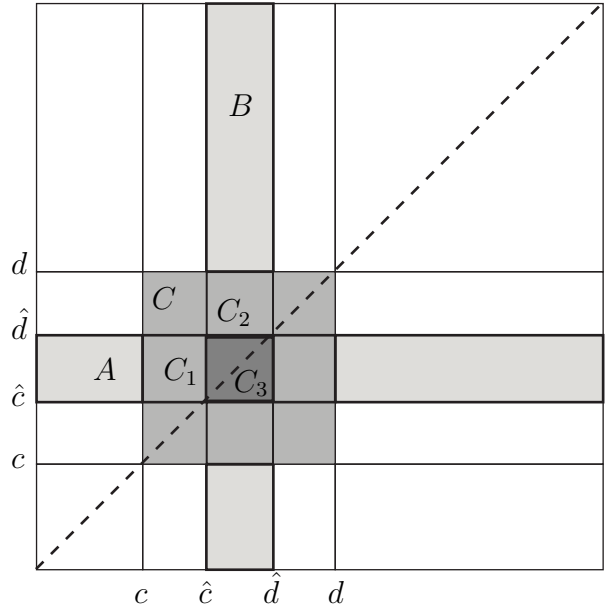}
                       \put(10,35){$A$}
                       \put(33,80){$B$}
                       \put(20,47){$C$}
                       \put(22, 35){$C_1$}
                       \put(31, 45){$C_2$}
                       \put(34,33){$C_3$} 
                       \put(17,-3){$c$}
                       \put(28,-3){$\hat{c}$}
                       \put(39,-3){$\hat{d}$}
                       \put(50,-3){$d$}
                       \put(-4,20){$c$}
                       \put(-4,30.5){$\hat{c}$}
                       \put(-4,42){$\hat{d}$}
                       \put(-4,53){$d$}
                       \end{overpic}}
    \caption{We break up our proof into several pieces according to a subdivision of the parameter space. Define $A:=[a,c]\times[\hat{c}, \hat{d}]$, $B:=[\hat{c},\hat{d}]\times[d,b]$, $C:=[c,d]\times[c,d]$, $C_1:=[c, \hat{c}]\times[\hat{c}, \hat{d}]$, $C_2:=[\hat{c}, \hat{d}]\times[\hat{d},d]$,  and $C_3:=[\hat{c}, \hat{d}]\times[\hat{c},\hat{d}]$. We use Lemma~\vref{DSbound} to control the change in the distortion quotient inside the region $C$. Then, we restrict our attention to the regions $A$ and $B$, which are bounded away from the diagonal of $[a,b]\times[a,b]$, and use a more straightforward calculation to show that the change in the distortion quotient can be controlled on regions $A$~and~$B$.}
    \label{reference}
  \end{figure}
  
\begin{case}
Assume the arc $\gamma((c,d))$ has total curvature $0<K~<~2 \cos^{-1}(1/(1~+~\varepsilon/2))$; then by Lemma~\vref{DSbound}, the distortion of the arc $\gamma((c,d))$ is no more than $1 + \varepsilon/2$.
\end{case}
\begin{proof}
By replacing the interval $\gamma([c,d])$ with an inscribed polygon, we can decrease length and not increase total curvature. By Lemma~\vref{DSbound}, since the total curvature has not increased, the distortion of the arc remains between 1 and $1 + \varepsilon/2$. Hence the change in distortion of the arc is smaller than $\varepsilon$. This takes care of all points in region $C$ in Figure~\ref{reference}.
  
Now let us turn our attention to a subarc $\gamma([\hat{c},\hat{d}]) \subset \gamma((c,d))$. We define $[\hat{c}, \hat{d}]$ to be any interval contained in $(c,d)$ so that $\gamma([\hat{c}, \hat{d}])$ has positive curvature. There exists a Radon measure $\mathcal{K}$ on $[a,b]$ so that $\mathcal{K}((p,q))$ is the total curvature of $\gamma((p,q))$ for all intervals $(p,q) \subset (a,b)$. We will use $\mathcal{K}$ to show that we can find an interval $[\hat{c}, \hat{d}]$ with the desired properties. Since $\mathcal{K}$ is a Radon measure, if we let $H$ denote the set of all compact subsets of $(c,d)$ then
\begin{equation*}
  \mathcal{K}((c,d)) = \sup_{h \in H} \mathcal{K}(h)
\end{equation*}
Hence there must exist some sequence of compact sets $h_i$ so that $\mathcal{K}(h_i) \rightarrow \mathcal{K}((c,d))$ as $i \rightarrow \infty$ and so there must also be some $j$ so that $\mathcal{K}(a_j) > 0$. Let $[\hat{c}, \hat{d}]\subset (c,d)$ be any closed interval containing $h_j$. Notice that we can replace $\gamma[\hat{c}, \hat{d}]$ with a polygon without increasing the distortion too much for points in $C_1$, $C_2$, and $C_3$ in Figure~\ref{reference} since each of these regions is contained within the larger region $C$.
 
For any given $\varepsilon' > 0$, let $P\colon[a,b]\longrightarrow \R^3$ be a curve that consists of two types of arcs. The first is any polygon inscribed inside the arc $\gamma([\hat{c},\hat{d}])$ such that $\| \gamma(s) - P(s) \| < \varepsilon'$ for all $s \in [\hat{c}, \hat{d}]$, $P(\hat{c}) = \gamma(\hat{c})$, and $P(\hat{d}) = \gamma(\hat{d})$. The second arc satisfies the relation $P(s) = \gamma(s)$ for all $s \in [a,b]\setminus [\hat{c}, \hat{d}]$. Even if the arc $\gamma([\hat{c}, \hat{d}])$ is already a polygon, we can replace $\gamma([\hat{c}, \hat{d}])$ with $P([\hat{c}, \hat{d}])$ so as to decrease length. 

It remains to show that we can choose $\varepsilon'$ small enough on $[\hat{c}, \hat{d}]$ so that 
\begin{equation*}
  dq_P(s,t) - dq_{\gamma}(s,t) < \varepsilon \text{ for all pairs } (s,t) \in A \cup B
\end{equation*}
where $A$ and $B$ are defined in Figure~\ref{reference}.

Notice that if $dq_P(s,t) < dq_{\gamma}(s,t)$ then there is nothing to show, since we are only concerned with \emph{increasing} distortion. So we will assume that $dq_P(s,t) > dq_{\gamma}(s,t)$. Suppose that $(s,t)~\in~A$. Then, using the definition of distortion quotient, 
\begin{equation*}
    dq_P(s,t) - dq_{\gamma}(s,t) =\frac{d(\gamma(s),\gamma(t);\R^3) \,d(P(s),P(t); P) - d(P(s), P(t); \R^3)\,d(\gamma(s), \gamma(t); \gamma)}{d(P(s), P(t); \R^3) \,d(\gamma(s),\gamma(t); \R^3) }
\end{equation*}

But $(s,t) \in A$ and since $A$ is bounded away from the diagonal of $[a,b] \times [a,b]$, we know that $s$ and $t$ must be at least $\min\{\hat{c} - c, d - \hat{d}\}>0$ apart. Thus, since the curves $\gamma$ and $P$ are both embeddings, we know 
there is a constant $C$ so that
\begin{equation*}
C \geq \frac{1}{d(P(s),P(t);\R^3)\,d(\gamma(s),\gamma(t);\R^3)}
\end{equation*}
for all $(s,t)\in A$. A calculation shows that
\begin{equation*}
     \frac{dq_P(s,t)  - dq_{\gamma}(s,t)}{C} \leq d(\gamma(s), \gamma(t); \gamma) \left[d(\gamma(s),\gamma(t);\R^3)  -  d(P(s),P(t); \R^3) \right].
 \end{equation*}
 But, $d(P(t), \gamma(s); \R^3) \leq  d(P(t), P(s); \R^3) + d(P(s), \gamma(s); \R^3)$ and by rearranging and rewriting terms $- d(P(s), P(t); \R^3) \leq d(P(s), \gamma(s); \R^3) - d(P(t), \gamma(s); \R^3)$. Using these facts, we see 
 \begin{equation*}
   d(\gamma(s),\gamma(t);\R^3)  -  d(P(s),P(t); \R^3) \leq d(P(t),\gamma(t); \R^3) + d(P(s), \gamma(s); \R^3) \leq 2\varepsilon'.
\end{equation*}

Finally, since $s,t \in [\hat{c}, \hat{d}]$, we can choose $\varepsilon'$ so that 
\begin{equation*}
  2\varepsilon' < \varepsilon/(Cd(\gamma(s),\gamma(t);\gamma)),
\end{equation*}
completing the proof in this case.
\end{proof}

\begin{case}Suppose there is no subarc of $\gamma((c,d))$ with total curvature $K$ where $$0~<~K~<~2\cos^{-1}(1/(~1~+~\varepsilon/2))$$.
\end{case}
\begin{proof} We will show that this forces $\gamma((c,d))$ to be a polygonal arc. A new technique will be required in this case.

By the Lebesgue-Radon-Nikodym theorem, we can write the curvature measure $\mathcal{K}$ as a sum of two measures $\lambda_a$ and $\lambda_s$ so that $\lambda_a$ is absolutely continuous with respect to $ds$ and $\lambda_s$ is singular with respect to $ds$. We then further decompose $\lambda_s$ as a sum of atomic measures $\{\mu_{x_i}\}$ and a non-atomic measure $\mu$ which is singular with respect to arclength. So we have
\begin{equation*}
  \mathcal{K}((p,q)) = \lambda_a((p,q)) + \sum_i \mu_{x_i}((p,q)) + \mu((p,q)).
\end{equation*}

Suppose now that $K_0:=2\cos^{-1}(1/(~1~+~\varepsilon/2))$ and $(p,q)$ is an interval so that $\mathcal{K}((p,q)) > K_0$. Then $\sum_i \mu_{x_i}((p,q)) < \infty$ and so there exists a constant $M$ so that $\sum_{i > M}\mu_{x_i}((p,q)) < K_0/2$. Without loss of generality, assume $x_1 < \ldots < x_M$. Then we partition the interval $(p,q)$ with the points $\{x_1, x_2, \ldots, x_M\}$ as seen in Figure~\ref{deltas_unordered_pic}.

 \begin{figure}[ht!]
    \centerline{\includegraphics[]{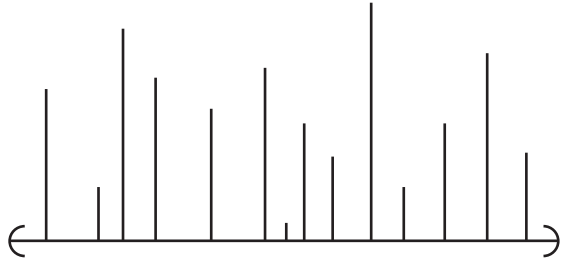}}
    \caption{Here we see a partition of the interval $(c,d)$ where the points that define the partition are the first $M$ atoms in the sequence $\{x_i\}$.}
    \label{deltas_unordered_pic}
  \end{figure}

We now make use of the following property of non-atomic Radon measures.

\begin{lemma}\label{measuregoodness}
  Suppose that $\mu$ is a non-atomic Radon measure (i.e., $\mu(\left\{p\right\}) = 0$ for each~$p\in\R$)  defined on $(a,b)$. Then given any finite interval $(a,b)$, we can find a value $L \in \R$ so that the measure of every subinterval $(c,d)$ of $(a,b)$ with length less than $L$ satisfies the inequality $\mu((c,d)) \leq 2\mu((a,b))/3$.
\end{lemma}

Now, let $(u,v)$ be a subinterval of $(p,q)$ that lies between two adjacent partition points. Since the measure $\lambda_a + \mu$ is non-atomic, by Lemma~\vref{measuregoodness} there is a value $L$ so that if $|v-u| < L$, then $\lambda_a((u,v)) + \mu((u,v)) < K_0/2$. Thus, we may assume that $(u,v)$ has small enough length so that  $\lambda_a((u,v)) + \mu((u,v)) < K_0/2$. Since $(u,v)$ does not include any of the partition points $x_1, \ldots, x_M$,
\begin{equation*}
\sum_i\mu_{x_i}((u,v)) = \sum_{i >M}\mu_{x_i}((u,v)) < K_0/2.
\end{equation*}
Thus $\mathcal{K}((u,v)) < K_0$ and it must be the case that $\mathcal{K}((u,v)) = 0$. Covering each interval between the $M$ partition points $\{x_1, x_2, \ldots, x_M\}$ by overlapping open intervals with length less than $|v-u|$ shows us that $\mathcal{K}((p,q) - \{x_1, \ldots, x_M\}) = 0$. Therefore, the curve is polygonal on the interval $(p,q)$ with corners at the atoms $\{x_1, x_2, \ldots, x_M\}$. Since this is true for any interval $(p,q)$ it must be the case that the arc $\gamma((c,d))$ is polygonal.

It remains to describe how to decrease length without changing the distortion quotient any more than $\varepsilon$.  First we will restrict our attention to a neighborhood of one corner small enough to guarantee that no pair of points in the neighborhood are $\Len(\gamma)/2$ apart along the curve $\gamma$. If the exterior angle at the corner point is $\phi$, then a calculation in \cite{mullikin} shows the restriction of the distortion quotient to the edges meeting at the corner achieves a maximum value of $\sec \frac{\phi}{2}$ as seen in Figure~\ref{polygon_case_new}.

  \begin{figure}[h]
    \centerline{\begin{overpic}{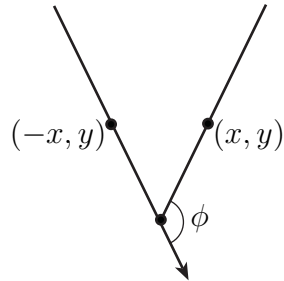}
                         \put(-15,50){$(-x,y)$}
                         \put(58,50){$(x,y)$}
                         \put(50,20){$\phi$}
                       \end{overpic}}
    \caption{If we translate the curve so that the vertex at the corner of the arc is symmetric about the $y$-axis, then a calculation shows that the distortion quotient restricted to this pair of segments is maximized at points of the form $(x,y)$ and $(-x,y)$. Moreover, this maximum is exactly equal to $\sec \frac{\phi}{2}$ where $\phi$ is the exterior angle at the corner. See \cite{mullikin} for a proof.}
    \label{polygon_case_new}
  \end{figure}
  
As seen in Figure~\ref{cornerdrc}, we replace the edge with vertices $(0,0)$ and $(x,y)$ with an edge with vertices $\epsilon(x,y)$ and $(x,y)$ (using a small value of $\epsilon$) and replace the edge with vertices $(-x,y)$ and $(0,0)$ with the edge with vertices $(-x,y)$ and $\epsilon(x,y)$. The triangle inequality guarantees that this alteration will decrease length. Since $\sec\left(\phi/2\right)$ is an increasing function on the interval $(0,\pi)$, decreasing the exterior angle decreases distortion.

  \begin{figure}[htbp]
    \centerline{\begin{overpic}{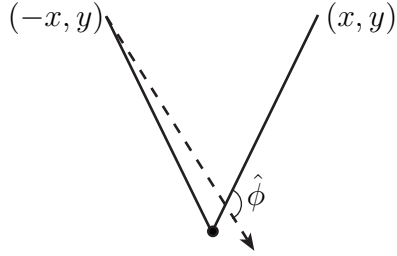}
                       \put(-40,95){$(-x,y)$}
                       \put(93,95){$(x,y)$}
                       \put(60,20){$\hat{\phi}$}
                       \end{overpic}}
    \caption{Here we shorten the length of the polygonal arc slightly while decreasing the exterior angle by a small amount by shortening one segment. This modification decreases the distortion quotient on this pair of line segments.}
    \label{cornerdrc}
    \end{figure}
This completes the proof.
\end{proof}
\end{proof}

\section{Main Theorem}
We now state the main theorem. Recall that distortion thickness $\tau_{\delta([\gamma])}(\gamma) \geq 1$ means that for any pair $(s,t)$ with $dq_{\gamma}(s,t) \geq \delta([\gamma])$ we have $\|\gamma(s) - \gamma(t)\| \geq 1$.
  \begin{thm} [Main Theorem] Let $U_C([\gamma])$ be the set of all finite total curvature curves $\gamma$ in $[\gamma]$, with distortion $\delta(\gamma) < C$ and distortion thickness $\tau_{\delta([\gamma])} \geq 1$. 

Then any open interval on a curve of minimum length in $U_C([\gamma])$ is either a straight line segment or contains an endpoint of a $\delta([\gamma])$-drc.
\end{thm}

\begin{proof}
  Let $\gamma \colon [a,b] \longrightarrow \R^3$ be a curve of minimum length in $U_C([K])$ and let $(c,d)$ be any open interval on $\gamma$. If the total curvature along $(c,d)$ is zero, then there is nothing to show. So assume that the total curvature along $(c,d)$ is positive. Further assume that the interval $(c,d)$ does \emph{not} contain an endpoint of a $\delta([\gamma])$-drc. We will deduce a contradiction. Since there is no endpoint of a $\delta([\gamma])$-drc contained within $(c,d)$, Proposition~\ref{prop2} states that there is a subinterval $(p,q)$ of $(c,d)$ so that if $(s,t) \in (p,q)\times[a,b]$, then $\delta([\gamma]) - dq_{\gamma}(s,t) > \varepsilon$ for some $\varepsilon > 0$. Provided the total curvature along $(p,q)$ is nonzero, Proposition \vref{prop1} enables us to decrease the length of the interval $(p,q)$ to obtain a new curve $\hat{\gamma}$ with the property that $dq_{\hat{\gamma}}(s,t) - dq_{\gamma}(s,t) < \varepsilon$ for $(s,t)~\in~(c,d)~\times~[a,b]~\cup~[a,b]~\times~(c,d)$. Thus no new $\delta([\gamma])$-drc's exist with endpoints in $(c,d)~\times~[a,b]~\cup~[a,b]~\times~(c,d)$. We claim that we have not increased the distortion quotient at any pair outside $(c,d)~\times~[a,b]~\cup~[a,b]~\times~(c,d)$. Indeed, if $(s,t)$ is any point outside $(c,d)~\times~[a,b]~\cup~[a,b]~\times~(c,d)$, then $d(\hat{\gamma}(s), \hat{\gamma}(t);\hat{\gamma})~\leq~d(\gamma(s), \gamma(t); \gamma)$ while $d(\hat{\gamma}(s), \hat{\gamma}(t);\R^3)~=~d(\gamma(s), \gamma(t); \R^3)$. Since $\delta(\gamma)~\geq~\delta([\gamma])$, the distortion of $\gamma$ is realized on such a pair and this implies that $\delta(\hat{\gamma})~\leq~\delta(\gamma)~<~C$. Therefore, $\hat{\gamma}~\in~U_C([\gamma])$ and $\Len(\hat{\gamma})~<~\Len(\gamma)$. This is the desired contradiction.
\end{proof}

\section{Conclusion and Future Directions}
The main theorem provides a good deal of structure for distortion minimizing curves provided they have a representative of shortest length. Indeed, the frequency of drcs is highly reminiscent of the results in \cite{criticality}. We hope this new information on drcs will aid in constructing better lower bounds for distortion. 

We point out that the main theorem is still somewhat unsatisfactory since we have not yet proved that there are  length-minimizing curves in $U_C([\gamma])$. Fixing this problem requires proving at least two conjectures:
\begin{conjecture}
For any curve $\gamma$ there exists a chord with distinct endpoints $\gamma(s)$ and $\gamma(t)$ so that $dq_{\gamma}(s,t) \geq \delta([\gamma])$. Hence we can rescale any knot in $[\gamma]$ to have $\delta([\gamma])$-distortion thickness 1.
\end{conjecture}

This conjecture means that the set $U_C([\gamma])$ will be easy to work with. In this case we further conjecture:

\begin{conjecture}
The set $U_C([\gamma])$ contains a minimizer for length for any knot type $[\gamma]$ and any $C > \delta([\gamma])$. 
\end{conjecture}

As pointed out earlier, we have a path towards an answer to Gromov's question. It remains to find an appropriate topological invariant and lower bound for distortion in terms of this invariant. A potential candidate for the topological invariant is a quantity known as the hull number defined in \cite{2ndhull}. The $n^{\text{th}}$-hull of a curve $\gamma$ is the set of all points $x$ for which any plane passing through $x$ intersects $\gamma$ in at least $2n$ points (we count a tangential intersection as a double intersection). Then, the hull number for a curve $\gamma$ is the largest value of $n$ for which the $n^{\text{th}}$-hull is nonempty. We make this a knot invariant by defining the hull number of the knot type $[\gamma]$ to be the infimal hull number of all knots in $[\gamma]$. 

The hull number is of particular interest since it has been shown in \cite{Ismestiev} that the hull number does not increase for the family $\{\gamma, \gamma \# \gamma, \gamma \# \gamma \# \gamma, \ldots\}$ and yet it does increase for a family of knot types which seems to be a prime candidate for having no universal upper bound for distortion (the $(n, n-1)$-torus knots). The minimum distortions of these knots are therefore particularly interesting.

In \cite{mullikin2} we will provide approximate locally distortion minimizing curves in these knot types, computed using simulated annealing.
 
\newpage
\section{appendix}
\setcounter{lemma}{3}
\begin{lemma}
  Suppose that $\mu$ is a non-atomic Radon measure (i.e., $\mu(\left\{p\right\}) = 0$ for each~$p\in\R$) defined on $(a,b)$. Then given any finite interval $(a,b)$, we can find a value $L \in \R$ so that the measure of every subinterval $(c,d)$ of $(a,b)$ with length less than $L$ satisfies the inequality $\mu((c,d)) \leq 2\mu((a,b))/3$.
\end{lemma}
\begin{proof}

  We will proceed by the method of contradiction. We assume that for all $L > 0$ there exists an interval $(c,d) \subset (a,b)$ with $d-c < L$ so that $\mu((c,d)) > 2\mu((a,b))/3$. For each positive integer $n$, let $L(n) = (b-a)/2^n$, and let $s_n$ be an open interval with length less than $L(n)$ so that  $s_n \subset (a,b)$ and $\mu(s_n) > 2\mu((a,b))/3$. Notice that if we define $S_n$ to be the closure of $s_n$, then $\mu(S_n) = \mu(s_n)$ since $\mu$ contains no atoms. Let $C_n = \{C_{n_i}\}$ be the closed cover of $[a,b]$ consisting of the sets $C_{n_i} = [a + (i-1)L(n), a + iL(n)]$ for $i = 1, 2, \ldots, 2^n$. So $C_1$ consists of two sets of equal length $C_{1_1} = [a, (a+b)/2]$ and $C_{1_2} = [(a+b)/2, b]$.
\begin{claim} \label{sublemma} There must be some integer $i$ so that $S_i \subset [a,(a+b)/2]$ or $S_i \subset [(a+b)/2, b]$.
\end{claim} 
\begin{proof}
  If there is no such integer, then it follows that $S_i \cap [a,(a+b)/2] \ne \emptyset$ and $S_i \cap [(a+b)/2, b] \ne \emptyset$ for all $i$. Therefore, $(a+b)/2 \in S_i$ for all $i$. Since the length of the $S_i$'s is approaching zero it follows that
\begin{equation*}
  \left\{\frac{(a+b)}{2}\right\} = \bigcap_{i = 1}^{\infty} S_i.
\end{equation*}
But then since $\mu$ is a Radon measure $\mu(\{(a+b)/2\}) > 2\mu((a,b))/3$, which contradicts the assumption that there are no atoms for $\mu$.
\end{proof}

Therefore, by Claim~\ref{sublemma} there exists some $S_i$ so that either $S_i \subset C_{1_1}$ or $S_i \subset C_{1_2}$. Hence, either $\mu(C_{1_1})>2\mu((a,b))/3$ or $\mu(C_{1_2})>2\mu((a,b))/3$. Define that set to be $T_1$. The next open cover $C_2$, covers $T_1$ with exactly two sets. Using the claim again and the fact that $\mu((a,b)\setminus T_1)~<~\mu((a,b))/3$, we know that only one member of $C_2$ has $\mu$-measure greater than $2\mu((a,b))/3$. Call it $T_2$. Continue inductively to generate a sequence of closed intervals $\{T_n\}$ with the following properties:

\begin{itemize}
  \item[(1)] $\mu(T_n) > 2\mu((a,b))/3$ for all $n$,
  \item[(2)] the length of $T_n$ is exactly $(b-a)/2^n$ for all $n$, and
  \item[(3)] $T_{n+1} \subset T_n$ for all $n \geq 1$, hence $\cap_{n=1}^{\infty} T_n \ne \emptyset$. 
\end{itemize}

Since the intersection of the $T_i$'s is nonempty and since the length of the $T_i$'s is approaching zero, it follows that the intersection of the $T_i$'s is a single point and that the $\mu$-measure of this point is larger than  $2\mu((a,b))/3 \geq 0$ as in Claim~\ref{sublemma}. This again contradicts the assumption that $\mu$ contains no atoms.
\end{proof}

\bibliographystyle{plain}
\bibliography{distortion}

\begin{thebibliography}{10}

\bibitem{criticality}
Jason Cantarella, Joseph~H.G. Fu, Robert~B. Kusner, Hohn~M. Sullivan, and
  Nancy~C. Wrinkle.
\newblock Criticality for the {G}ehring link problem.
\newblock arXiv:math.DG/0402212.

\bibitem{2ndhull}
Jason Cantarella, Greg Kuperberg, Robert~B. Kusner, and John~M. Sullivan.
\newblock The second hull of a knotted curve.
\newblock {\em Amer. J. Math}, 125(6):1335--1348, 2003.

\bibitem{DS}
Elizabeth Denne and J.M. Sullivan.
\newblock The distortion of a knotted curve.
\newblock arXiv:math.GT/0409438.

\bibitem{gromov}
Mikhael Gromov.
\newblock Filling {R}iemannian manifolds.
\newblock {\em J. Differential Geom.}, 18(1):1--147, 1983.

\bibitem{Ismestiev}
Ivan Ismestiev.
\newblock Hull number of torus knots.
\newblock arXiv:math.GT/0412139.

\bibitem{KS}
Robert~B. Kusner and John~M. Sullivan.
\newblock On distortion and thickness of knots.
\newblock In {\em Topology and geometry in polymer science (Minneapolis, MN,
  1996)}, volume 103 of {\em IMA Vol. Math. Appl.}, pages 67--78. Springer, New
  York, 1998.

\bibitem{milnor}
J.~W. Milnor.
\newblock On the total curvature of knots.
\newblock {\em Ann. of Math. (2)}, 52:248--257, 1950.

\bibitem{mullikin2}
Chad~A.S Mullikin.
\newblock Numerical approximation of local minima for distortion of knots using
  simulated ammealing.
\newblock in preparation.

\bibitem{mullikin}
Chad~A.S. Mullikin.
\newblock {\em On Length Minimizing Curves with Distortion Thickness Bounded
  Below and Distortion Bounded Above}.
\newblock PhD thesis, University of Georgia, 2006.

\bibitem{ohara2}
Jun O'Hara.
\newblock Family of energy functionals of knots.
\newblock {\em Topology Appl.}, 48(2):147--161, 1992.

\bibitem{Sullivan-FTC}
John~M. Sullivan.
\newblock Curves of finite total curvature.
\newblock arXiv:math.GT/0606007 v1, preprint.

\end{thebibliography}
\end{document}